\documentclass[11pt,francais]{smfart}
\usepackage{xspace,smfenum}
\usepackage[latin1]{inputenc}
\usepackage[T1]{fontenc}
\usepackage[frenchb]{babel}
\usepackage[cm]{aeguill}

\usepackage{mathrsfs}
  \let\mathcal\mathscr
\usepackage[matrix,arrow,curve,frame,cmtip]{xy}

\usepackage{smfthm}

\makeatletter
  \NumberTheoremsAs{paragraph}
  
  \numberwithin{equation}{paragraph}
  \@addtoreset{paragraph}{section}
  \@addtoreset{paragraph}{subsection}
  \def\theequation{\arabic{section}.\arabic{equation}}

  \def\subsection{\@startsection{subsection}{2}%
  \z@{.5\linespacing\@plus.7\linespacing}{.5\linespacing}%
  {\normalfont\bfseries}}

\makeatother

\def\Aut{\operatorname{Aut}}
\def\Aff/#1{(\mathrm{Aff}/#1)}
\def\id{\mathrm{id}}
\def\ob{\operatorname{ob}\nolimits}
\def\div{\operatorname{div}\nolimits}
\def\Fib{\text{Fib}}
\def\Isom{\operatorname{\mathbf{Isom}}}
\def\Spec{\operatorname{Spec}}
\def\Hom{\operatorname{Hom}}
\let\ra\rightarrow

\let\phi\varphi
\def\P{{\mathbf P}}
\def\C{{\mathbf C}}
\def\Q{{\mathbf Q}}
\def\V{{\mathbf V}}
\def\A{{\mathbf A}}
\def\propS #1{\textup{(S$_{#1}$)}}
\def\loccit{\emph{loc.\ cit.}\nobreak\xspace}
\def\cf{\emph{cf.}\nobreak\xspace}
\def\resp{\emph{resp.}\nobreak\xspace}

\begin{document}
\author{Antoine Chambert-Loir}
      \title{Champs de Hurwitz}
\address{CMAT, \'Ecole polytechnique \\ 91128 Palaiseau Cedex}
\email{chambert@math.polytechnique.fr}

\begin{abstract}
On construit les champs de Hurwitz et on en donne quelques
propriétés, essentiellement contenues dans SGA\,1.
Quelques applications de nature arithmétique en sont déduites.
\end{abstract}

\begin{altabstract}
We propose a construction of Hurwitz stacks and give some properties
of them, most of which are consequences of SGA\,1. Some applications
of arithmetic flavor are then deduced.
\end{altabstract}
\maketitle

\tableofcontents

\section{Introduction}

Les \emph{espaces de Hurwitz} sont, classiquement,
des variétés analytiques complexes qui paramètrent, en un sens
plus ou moins précis, les revêtements de la sphère de Riemann
dont le degré et le nombre de points de ramification sont donnés,
et dont la ramification est simple.
Ils ont notamment permis d'établir de manière algébrique
(Klein, Hurwitz) la connexité de l'espace
de modules des courbes en caractéristique zéro, et
même (Fulton, \cite{fulton69} en caractéristique positive, assez grande
(voir aussi l'appendice de~\cite{harris-m82} qui fournit
une démonstration algébrique du seul argument d'origine analytique
de la démonstration de Deligne et Mumford).
Leur construction était originellement par voie transcendante.
En particulier, l'article~\cite{emsalem95} d'Emsalem construit
effectivement ces espaces de Hurwitz en tant que variétés
algébriques sur le corps des rationnels, mais la démarche
qui y est suivie consiste à algébriser l'espace de Hurwitz
défini de manière transcendante, puis à descendre son corps
de définition de $\C$ à $\Q$, via la théorie de Weil.
Cependant, dans l'article déjà cité de Fulton, ils sont construits de manière
algébrique, en recourrant à un théorème de Grothendieck  concernant
la \emph{représentabilité des foncteurs non ramifiés}

L'étude du problème inverse de Galois (construire des extensions
régulières de $\Q(T)$ de groupe de Galois donné)
les a remis au goût du jour il y a une vingtaine d'années (travaux
de Fried, \cite{fried77}) avec nécessairement un parfum arithmétique
un peu plus poussé. Les notions de \emph{corps de modules}
ou de \emph{corps de définition} d'une courbe algébrique
(\resp d'un revêtement) étaient au centre d'un grand nombre
de travaux. Outre bien sûr ceux concernant la rigidité,
ou l'effervescence autour des dessins d'enfants,
je pense surtout à un article de Beckmann (\cite{beckmann89})
dont il sera question dans ce texte : elle y contrôle
les places de ramification dans le corps des modules d'un revêtement
ramifié de la droite  projective. En revanche, il est bien connu
qu'il est en général plus difficile de préciser un corps de définition.

\medskip

Cette note tire en effet son origine (en 1995) de la constatation
que le théorème principal de cet article \emph{était}
une conséquence relativement directe de la théorie de la spécialisation
du groupe fondamental, telle qu'étudiée par Grothendieck 
et Mme~Raynaud dans SGA~1 (\cite{sga1}, exposé XIII).
Suite à des discussions avec Emsalem, cette approche a été effectivement
décrite dans un article de ce dernier (\cite{emsalem1999}).
Comme il est clair à quiconque a lu le dernier exposé de SGA~1,
cette approche s'étend aux revêtements d'une variété
algébrique lisse ramifiés le long d'un diviseur à croisements normaux
stricts.

Peu après, j'ai réalisé qu'on pouvait en fait interpréter
ce théorème comme un analogue du théorème de Chevalley-Weil
pour ces espaces de Hurwitz.
Tout le sel de l'histoire vient alors du parallèlisme
entre corps de définition/corps des modules
et variété/champ algébrique.

Parallèlement à l'élaboration de cet article,
l'importance des gerbes et de la cohomologie non abélienne
dans ces questions a été remarquablement mise en valeur
dans des travaux de Dèbes, Douai et Emsalem
(\cite{debes-d1997}, \cite{debes-d1998}, \cite{debes-d1999},
\cite{debes-e1999}, \cite{debes-d-e2000}). 

Dans le présent article, je construis des \emph{champs de Hurwitz}
dans une situation très générale. Leur existence en tant
que champs algébriques d'Artin est assez facile. Le fait
que ce soient des champs de Deligne-Mumford provient essentiellement
des propriétés de la spécialisation du groupe fondamental.
Les espaces grossiers qui leur sont associés sont a priori
des espaces algébriques, mais, miraculeusement, 
la finitude du groupe fondamental en fait des schémas.
Faute d'avoir réalisé ce dernier point plus tôt (\string!),
j'avais laissé cet article hiberner, certes un peu trop longtemps.
Je montre pour finir comment retrouver effectiment le théorème
de Beckmann mentionné plus haut dans le style Chevalley-Weil.

Notons aussi 
que les gerbes construites dans les articles de Dèbes, Douai, Emsalem,
apparaissent alors (évidemment) comme les \emph{gerbes résiduelles}
du champ de Hurwitz en ses points algébriques.

\medskip

Je ne peux terminer cette introduction sans mentionner
un certain nombre de lacunes (sérieuses) de cette note :
\begin{itemize}
\item l'article~\cite{harris-m82} de Harris et Mumford, repris
amplement par Mochizuki (\cite{mochizuki95}) introduit une compactification
des espaces de Hurwitz. Ce dernier incorpore aussi
des techniques logarithmiques. 
La thèse (apparemment non publiée) de Wewers propose
aussi une construction (différente) des espaces de Hurwitz usuels.

Il n'en sera hélas pas question ici. On peut espérer que la théorie
de la spécialisation du groupe fondamental en géométrie
logarithmique récemment étudiée par Vidal,
jointe aux méthodes de cette note permette de construire
des {\og log-champs\fg} qui ne sont sûrement pas sans rapport avec
la compactification de Harris-Mumford.

\item 
Je ne tiens pas compte des groupes d'inertie.

\item la vertu essentielle des espaces de Hurwitz classiques
est leur connexité (cf.~\cite{fulton69})
ainsi qu'une note récente de Graber, Harris et Starr, \cite{graber-h-s2002}).
En général, on ne peut probablement rien dire, mais je n'ai
pas réfléchi au problème.

\item en certaines caractéristiques, les champs de Hurwitz
que je construis ne sont probablement pas raisonnables, et d'ailleurs,
j'ignore largement ce qui se passe en ces caractéristiques. (C'est lié au fait
qu'il y a souvent trop de $p$-revêtements en caractéristique~$p$.)
Par d'autres méthodes (espaces de modules des courbes stables),
Abramovich et Oort ont défini dans~\cite{abramovich-o2000}
des compactifications des
espaces de Hurwitz.
\end{itemize}

\bigskip

Je remercie Dèbes, Emsalem, Lochak et Pop de l'intérêt
qu'ils ont porté à ce travail.

\section{Notations}

\paragraph{}
Nous utiliserons dans ce texte la théorie des champs algébriques,
telle qu'elle est développée dans la monographie~\cite{laumon-mb2000}
de Laumon et Moret-Bailly.
Depuis la rédaction de cet ouvrage, on peut noter que le critère
d'Artin d'algébricité d'un $S$-champ (corollaire~10.11 de \loccit)
a été étendu par Conrad et de Jong
du cas d'un schéma localement de type fini sur un corps ou sur un anneau de
Dedekind excellent au cas d'un schéma excellent quelconque.
Même si cette généralité est illusoire dans les applications
que nous avons en vue, nous rédigerons cette note en tenant compte
de cette amélioration.

\paragraph{}
Soit $S$ un schéma excellent.

On désigne par $\Aff/S$ la catégorie des schémas affines $T$
munis d'un morphisme $T\ra S$ ; cette catégorie est munie
de la topologie {\itshape fppf}. Conformément aux conventions
choisies dans~\cite{laumon-mb2000}, nous pourrions nous contenter
de la topologie étale. D'après \loccit, (9.4), ça ne change rien.

Si $X$ et $T$ sont deux $S$-schémas, on désigne
par $X_{T}$ le $T$-schéma $X\times_S T$ (changement de base).

Si $Y\ra X$ est un morphisme de schémas, on désigne par $\Aut(Y/X)$
l'ensemble des $X$-automorphismes de $Y$.

Un \emph{revêtement} d'un schéma $X$ est
un morphisme $Y\ra X$ localement libre de type fini et surjectif.

\section{Construction des champs de Hurwitz}

\paragraph{}
Soient $S$ un sch\'ema  excellent, ainsi qu'un sch\'ema $X$ muni d'un morphisme
projectif et plat $\pi:X\ra S$ à fibres géométriques connexes
et réduites.
(Donc $\pi_*\mathcal O_X=\mathcal O_S$, universellement.)

On s'intéresse à paramétrer les rev\^etements (localement libres
de type fini, surjectifs) de $X$,
vérifiant éventuellement certaines propriétés supplémentaires
(ramification, degré, groupe d'automorphisme\dots).

\begin{defi}
On définit une catégorie $\mathcal H(X/S)$ comme suit :
\begin{itemize}
\item
ses objets sont les triplets $(T,Y,f)$,
où $T\in\ob\Aff/S$, $Y$ est un schéma et $f:Y\ra X_{T}$ est
un revêtement ;
\item
un morphisme de $(T,Y,f)$ dans $(T',Y',f')$ est un couple
$(h,\eta)$, où $h:T\ra T'$, $\eta:Y\ra Y'$ sont deux morphismes
de schémas tels que le diagramme
\[ \xymatrix{
     Y \ar^{\eta}[r] \ar[d] & Y' \ar[d] \\
     X_T \ar^{\id_X\times h}[r] & X_{T'} } \]
est cartésien.
\end{itemize}
On munit cette catégorie du foncteur d'{\og oubli du rev\^etement\fg}
$\mathcal H(X/S)\ra \Aff/S$
qui associe à $(T,Y,f)$ le schéma $T$.
\end{defi}

Quand il ne pourra y avoir de confusion possible, on notera
$\mathcal H=\mathcal H(X/S)$. Si $T$ est un $S$-schéma,
on notera aussi $\mathcal H_T$ la catégorie fibre en~$T$
(dont les objets sont les $(T,Y,f)$ comme ci-dessus).

\begin{lemm}
La catégorie $\mathcal H(X/S)$ ainsi définie est un $S$-champ.
\end{lemm}
\begin{proof}
Il est immédiat sur la définition que $\mathcal H$
est un $S$-groupo{\"\i}de (\cite{laumon-mb2000}, définition 2.1).
Il y a alors deux points \`a v\'erifier~:
\begin{itemize}
\item
Soient $T\ra S$ un sch\'ema, $f_1:Y_1\ra X_T$ et $f_2:Y_2\ra
X_T$ deux objets de $\mathcal H_T$. Alors, le foncteur
$\Isom_{T}(Y_1,Y_2)$ est un faisceau sur $T$. C'est clair puisque
premi\`erement, les isomorphismes dans $\mathcal H_T$ sont
des morphismes de $T$-sch\'emas et qu'un morphisme entre deux $T$-sch\'emas
se d\'efinit localement sur $T$~; deuxi\`emement, un tel morphisme
qui est localement cart\'esien sur $T$ l'est globalement.
\item
Si $T'\ra T$ est un morphisme fppf, toute donn\'ee de descente de $T'$
\`a $T$ dans $\mathcal H$ est effective.
Soit en effet $f':Y'\ra X\times_T {T'}$
un objet de $\mathcal H_{T'}$, muni d'une donn\'ee de descente de $T'$ \`a
$T$. Un morphisme fini \'etant
affine, il r\'esulte du th\'eor\`eme de descente de Grothendieck
(\cite{sga1}, exp.~\textsc{\romannumeral 8}, th.~2.1)
qu'il existe un objet canonique $f:Y\ra X_T$ de $\mathcal H_{T}$
tel que $f\times_{T}T'$ soit isomorphe \`a $f'$.
D'autre part, la propri\'et\'e pour un morphisme d'\^etre 
fini localement libre est locale pour la topologie fppf
(\loccit, cor.~5.7).
\qed
\end{itemize}
\let\qed\relax
\end{proof}

\paragraph{}\label{parag.varH}
Introduisons maintenant diverses généralisations du champ $\mathcal H(X/S)$.
\begin{itemize}
\item
Soit $n$ un entier $\geq 1$ ; la catégorie $\mathcal H^n(X/S)$ est
la sous-catégorie pleine de $\mathcal H$  formée des
$(T,Y,f:Y\ra X_T)\in\mathcal H_T$ tels que $f$ soit de degré constant $n$.
\item
Soit $U$ un ouvert de $X$ schématiquement
dense dans toute fibre de $X\ra S$, si bien que
pour tout sch\'ema $T\ra S$, $U\times_S T$ est
schématiquement dense dans $X\times_S T$ ; la catégorie $\mathcal H^{\text{$U$-ét}}(X/S)$
est la sous-catégorie pleine de $\mathcal H$ formée
des $(T,Y,f)\in\mathcal H_T$ tels que $f$ soit étale au-dessus
de $U\times_S T$.
\item
La catégorie $\mathcal H^{\text{g-ét}}(X/S)$ est la sous-catégorie pleine
de $\mathcal H$ formée des $(T,Y,f)\in\mathcal H_T$ tels que
l'ensemble des points de $X_T$ au-dessus desquels 
$f$ est étale soit un ouvert de $X_T$ qui est schématiquement dense dans
toute fibre de $X_T\ra T$.
\end{itemize}

On peut aussi ajoindre des structures supplémentaires aux
rev\^etements, comme par exemple :
\begin{itemize}
\item
Soit $G$ un groupe fini ; la catégorie $\mathcal H^{G}(X/S)$ a pour objets
les $((T,Y,f),\phi)$, où $(T,Y,f)$ est un objet de $\mathcal H$
et $\phi:G\ra\Aut(Y/X_T)^{\text{op}}$ est un anti-isomorphisme
de $G$ dans les automorphismes de $Y/X_T$. Les morphismes de la catégorie
$\mathcal H^G$ sont les morphismes de la catégorie $\mathcal H$
qui commutent à l'action de $\phi$.
\item
Si $Z\subset X$ est un sous-schéma tel que $\pi_*\mathcal O_Z=\mathcal O_S$
universellement, et si $n$ est un entier $\geq 1$,
la catégorie $\mathcal H^{n,\text{$Z$-tr}}(X/S)$ a pour objets
les $((T,Y,f),\sigma)$, où $(T,Y,f)$ est un objet
de $\mathcal H$ et $\sigma$ est un isomorphisme
$\{1,\ldots,n\}\times Z_T \ra f^ {-1}(Z_T)$.
\end{itemize}

On définit d'autres catégories $\mathcal H^{n,\text{$U$-ét}}$,
$\mathcal H^{G,\text{$U$-ét}}$, $\mathcal H^{n,\text{$U$-ét},\text{$Z$-tr}}$,
etc.\ qui consistent à imposer toutes
les conditions et à adjoindre toutes les données indiquées en exposant,
les morphismes étant bien s\^ur astreints à vérifier les
compatibilités correspondantes.
En outre, on dispose de foncteurs d'oubli naturels vers $\mathcal H(X/S)$ qui
permettent d'identifier les catégories avec plusieurs conditions comme
les produits fibrés des catégories correspondantes
au-dessus de $\mathcal H$.

Dans les applications que nous avons en vue, {\em $X$ sera en outre lisse
sur $S$ et $U$ sera le compl\'ementaire d'un diviseur \`a croisements
normaux relatif sur $S$.}
Le cas le plus simple, donnant lieu aux espaces de
Hurwitz classiquement consid\'er\'es, est obtenu pour
$S$ l'espace qui paramètre les familles de $d$ points distincts
de la droite projective sur~$\Q$, $X=\P^1\times S$,
$D\subset X$ le diviseur de degré~$d$ universel
(dont la restriction à $\P^1\times\{s\}$ \emph{est} le diviseur
donné par le point $s$) et $U=X\setminus D$ l'ouvert complémentaire.

\begin{lemm}\label{aut.affine}
Pour tout objet $(T,T,f:Y\ra X_T)$ de $\mathcal H(X/S)$, le faisceau
sur $\Aff/T$ donné par $T'\mapsto \Aut(Y_{T'}/X_{T'})$
est représentable par un $T$-schéma affine de type fini.
\end{lemm}
\begin{proof}
Soit $\mathcal E=f_*\mathcal O_Y$ la $\mathcal O_{X_T}$-algèbre localement
libre de rang fini
définissant $Y$ (de sorte que $Y=\Spec_{X_T} \mathcal E$). Il y a
bijection entre $\Aut(Y_{T'}/X_{T'})$ et les automorphismes
de $\mathcal E_{T'}$ comme $\mathcal O_{X_{T'}}$-algèbre, d'où
un faisceau représentable par un sous-schéma fermé de
$ \V(\pi_*\Hom_{\mathcal O_{X_T}}(\mathcal E,\mathcal E))^2$.
\end{proof}


\begin{prop}
Toutes ces catégories introduites sont des $S$-champs.
De plus, les $1$-morphismes naturels
\[ \mathcal H^n(X/S)\ra\mathcal H(X/S),
\quad 
   \mathcal H^{\text{$U$-ét}}(X/S)\ra\mathcal H^{\text{g-ét}}(X/S)
   \ra\mathcal H(X/S) \]
sont représentables par des immersions ouvertes.
Les $1$-morphismes
\[ \mathcal H^{G}(X/S)\ra \mathcal H(X/S)
\quad\text{(resp.\ } \mathcal H^{\text{$Z$-triv}}(X/S)\ra\mathcal H(X/S)
\text{)} \]
sont représentable affines, et le second est m\^eme
un $\mathfrak S_n$-torseur.
\end{prop}
\begin{proof}
Une sous-catégorie pleine d'un $S$-groupo\"\i de en est automatiquement un,
si bien que $\mathcal H^n$, $\mathcal H^{\text{$U$-ét}}$ et
$\mathcal H^{\text{g-ét}}$ sont des $S$-groupo\"\i des.
De m\^eme, pour ces $S$-groupo\"\i des, la propriété que
le foncteur des isomorphismes entre deux objets soit un faisceau
résulte de la propriété correspondante dans $\mathcal H(X/S)$.

D'autre part, les conditions imposées sont locales pour la topologie
fppf,
si bien que toute donnée de descente qui est déjà effective
dans $\mathcal H(X/S)$
fournira un objet de la sous-catégorie, laquelle est ainsi un
sous-$S$-champ.

Pour prouver que le $1$-morphisme $\mathcal H^n\ra\mathcal H$
est représentable
par une immersion ouverte, il faut prouver que pour tout
$S$-schéma $T$, et tout $Y\ra X_T$, l'ensemble des
$t\in T$ tels que $Y_t\ra X_t$ est de degré $n$ est ouvert dans $T$.
Or, le degré d'un morphisme localement libre est localement constant
si bien que $Y\ra X_T$ est de degré $n$ au-dessus d'une
réunion de composantes connexes de $X_T$ lesquelles sont
de la forme $X_{T_1}$, pour $T_1$ une composante connexe de $T$,
puisque $X\ra S$ est à fibres géométriquement connexes.

De m\^eme, soit $f:Y\ra X_T$ un rev\^etement.
Soit $B\subset X_T$ son lieu de branchement, plus petit fermé
de $X_T$ tel que $f:Y\ra X_T$ soit étale au-dessus de $X_T\setminus B$.
Comme $\pi:X\ra S$ est plat, la projection $\pi(X_T\setminus B)$ 
est un ouvert $T'\subset T$.
C'est de plus le plus grand sous-schéma de $T$ tel
que $f:Y_{T'}\ra X_{T'}$ est étale au-dessus d'un ouvert de $X_{T'}$
surjectif sur $T'$, et sa formation commute au changement de base.
Par conséquent, $\mathcal H^{\text{g-ét}}\ra\mathcal H$ est représentable
par une immersion ouverte.

Finalement, soit $U\subset X$ un ouvert de $X$ 
qui est dense dans toute fibre de $X\ra S$ 
et montrons que $\mathcal H^{\text{$U$-ét}}\ra\mathcal H$
est représentable par une immersion ouverte.
Considérons donc $f:Y\ra X_T$ un rev\^etement.
Son lieu de branchement est un fermé $B\subset X_T$ ;
l'intersection $B\cap (X_T\setminus U)$ est un fermé de $X_T$
dont la projection dans $T$ est fermée. Notant $T'$ l'ouvert
complémentaire, on constate que $f_{T'}:Y_{T'}\ra X_{T'}$ est
un objet de $\mathcal H^{\text{$U$-ét}}$ et $T'$ est le plus grand
sous-schéma de $T$ qui vérifie cette propriété, d'où l'assertion.

\medskip

Passons maintenant au champ $\mathcal H^G$.
C'est tout d'abord bien un $S$-champ : il est déjà clair que c'est
un $S$-groupo\"\i de. De plus le préfaisceau des isomorphismes
entre deux objets est un faisceau : étant donnés deux
objets $(Y_1\ra X_T;\phi_1)$ et $(Y_2\ra X_T;\phi_2)$
de $\mathcal H^G_T$ et un morphisme fppf $T'\ra T$, tout isomorphisme
$\alpha : (Y_{1,T'}\ra X_{T'};\phi_{1,T'})\ra (Y_{2,T'}\ra
X_{T'};\phi_{2,T'})$ dans $\mathcal H^G_{T'}$ descend en un isomorphisme
dans $\mathcal H_{T'}$ qui va nécessairement commuter aux actions
de $G$ via  $\phi_1$ et $\phi_2$.
Finalement, montrons que toute donnée de descente est effective :
une donnée de descente de $T'$ à $T$ sur un objet
$(Y_{T'}\ra X_{T'};\phi_{T'})$ de $\mathcal H^G_{T'}$
nous fournit tout d'abord un objet $Y\ra X_T$ de $\mathcal H_T$
car $\mathcal H_T$ est un $S$-champ.
On obtient une donnée de descente de $T'$ à $T$ sur les $\phi_{T'}(g)$
qui est elle aussi effective, $\Aut(Y/X_T)$ étant un $T$-schéma,
d'où un anti-homomorphisme $\phi:G\ra \Aut(Y/X_T)^{\text{op}}$
qui est, encore par descente, un anti-isomorphisme.

La représentabilité du $1$-morphisme $\mathcal H^G\ra\mathcal H$
résulte du fait que pour toute famille $(Y\ra X_T)\in\mathcal H_T$,
le foncteur des homomorphismes $T'\mapsto \Hom(G_{T'}\ra\Aut(Y_{T'}/X_{T'})$ 
est le foncteur des $T$-homomorphismes entre deux
schémas affines sur $T$,
et est donc représentable par un schéma affine affine sur~$T$.

\medskip
Enfin, montrons que $\mathcal H^{n,\text{$Z$-tr}}$ est
un champ et que $\mathcal H^{n,\text{$Z$-tr}}\ra\mathcal H$
est représentable, fini et plat.
Qu'il soit un $S$-champ se prouve de manière analogue (mais plus
simple) à $\mathcal H^G$.
Quant à la représentabilité, si $Y\ra X_T$ est un objet de $\mathcal H_T$,
soit $\mathcal E$ la $\mathcal O_{X_T}$-algèbre localement libre
correspondant à $Y$.
Les $\sigma:\{1,\ldots,n\}\times Z_T\ra f^ {-1}(Z_T)$ s'identifient
à des isomorphismes de $\mathcal O_{Z_T}$-algèbres
\[ \mathcal E|_{Z_T} \simeq \mathcal O_{Z_T}^{\oplus n} , \]
donnée représentable par un torseur sous le groupe
symétrique $\mathfrak S_n$, étant donné que $\pi_*\mathcal O_Z=\mathcal
O_S$ universellement.
\end{proof}

\begin{prop}\label{prop.champ.alg.tf}
Soit $n\geq 1$ un entier.
Le champ $\mathcal H^n(X/S)$ est un champ algébrique de type fini.
Par conséquent, le champ $\mathcal H(X/S)$
et ses variantes introduites en (\ref{parag.varH})
sont des champs algébriques localement de type fini.
\end{prop}
\begin{proof}
Soit $\Fib_n$ le $S$-champ des faisceaux localement libres
de rang~$n$ sur $X$. D'après le théorème~4.6.2.1
de~\cite{laumon-mb2000},
c'est un champ algébrique de type fini sur~$S$.
On dispose d'un $1$-morphisme de $S$-champs
$\mathcal H_n\ra\Fib_n$ qui associe
au rev\^etement de degré $n$, $f:Y\ra X_T$, le 
faisceau localement libre de rang~$n$ sur $X_T$ qu'est $f_*\mathcal O_Y$.

Nous alons prouver que ce morphisme est représentable.
Compte tenu du dictionnaire entre schémas affines sur $X_T$
et $\mathcal O_{X_T}$-algèbres quasi-cohérentes,
cela revient à prouver qu'étant donné un faisceau 
$\mathcal E$, localement libre de rang~$n$ sur $X_T$,
la catégorie formée des structures de $\mathcal O_{X_T}$-algèbres
sur $\mathcal E$ est représentable par un espace algébrique.
Or, elle l'est m\^eme par un schéma ! Une structure de $\mathcal O_{X_T}$
algèbre sur $\mathcal E$ consiste en effet en :
\begin{itemize}
\item
une unité,
qui est une section ne s'annulant pas de $\pi_*\mathcal E$, représentée
par l'ouvert complémentaire de la section nulle dans le $T$-fibré vectoriel
$\V(\pi_*\mathcal E^\vee)$ ;
\item
une fois fixée l'unité,
une loi de multiplication, astreinte aux compatibilités qui font
de $\mathcal E$ une $\mathcal O_{X_T}$-Algèbre, représentée par un
fermé du $T$-fibré vectoriel 
$\V(\pi_*\Hom(\mathcal E\otimes\mathcal E,\mathcal E)^\vee)$.
\end{itemize}
La réunion de ces deux conditions est ainsi représentable par
un schéma quasi-projectif de type fini, d'où la proposition.
\end{proof}

\begin{prop}
On suppose que les fibres du  morphisme $X\ra S$ satisfont à la propriété
\propS2 de Serre\footnote{En d'autres termes, le morphisme $X\ra S$
est \propS2, \cf EGA~4, déf.~6.8.1.}.
Soient $T$ un $S$-schéma, $x_1=(f_1:Y_1\ra X_T)$ et $x_2=(f_2:Y_2\ra X_T)$
deux objets de $\mathcal H^n_T$.
Supposons qu'il existe un ouvert $U\subset X_T$
schématiquement dense dans toute fibre au-dessus de~$T$ 
tel que $f_1$ et $f_2$ soient \'etales
au-dessus de $U$. Alors, le faisceau sur $T$
donné par $\Isom_T(x_1,x_2)$ est un schéma fini
et non-ramifié sur $T$.
\end{prop}
\begin{proof}
On peut bien sûr supposer que $T=S$.
La représentabilité de ce faisceau est contenue dans la définition d'un
champ algébrique, ainsi que le fait qu'il soit séparé
et quasi-compact sur $S$. Nous avons
déjà vu qu'il est affine (lemme~\ref{aut.affine}). Montrons
qu'il est propre et non ramifié.

Soient $Y_1\ra X$ et $Y_2\ra X$
les rev\^etements de degr\'e $n$ correspondant
\`a $x_1$ et $x_2$.

\emph{Prouvons que le faisceau des isomorphismes
entre $x_1$ et $x_2$ est non ramifi\'e.}
On peut pour cela supposer
que $S=\Spec R$, $I\subset R$ est un id\'eal de $R$ tel que $I^2=0$
et noter $S_0=\Spec(R/I)$.
\'Etant donn\'e un isomorphisme $\alpha_0:Y_{1,S_0}\ra Y_{2,S_0}$
au-dessus de $X_{S_0}$, il faut montrer qu'il existe au plus
un isomorphisme $\alpha:Y_1\ra Y_2$ au-dessus de $X$
dont $\alpha_0$ soit la restriction \`a $S_0$.

La {\og remarquable \'equivalence de cat\'egories\fg}
(EGA 4, th.~18.1.2)
implique que la restriction de $\alpha$ \`a $f_1^{-1}(U)$
est uniquement d\'etermin\'ee par la restriction de $\alpha_0$
\`a $f_1^{-1}(U)_{S_0}$. Comme $U$ est schématiquement dense dans $X$,
il existe
au plus un isomorphisme $\alpha:Y_1\ra Y_2$ qui rel\`eve $\alpha_0$.

\emph{Prouvons qu'il est propre.}
On peut alors supposer que $S$ est un trait de point
g\'en\'erique $\eta$, de point ferm\'e $s$. Si $\alpha_\eta:Y_{1,\eta}
\ra Y_{2,\eta}$ est un isomorphisme au-dessus de $X_\eta$,
montrons qu'il existe un unique isomorphisme $\alpha:Y_1\ra Y_2$
dont la restriction \`a $\eta$ soit $\alpha_\eta$.
On peut interpréter $\alpha_\eta$ comme une section du morphisme
canonique $Y_1\times_X Y_2 \ra Y_1$ au-dessus de $Y_{1,\eta}$
qu'il faut étendre à $Y_1$.

D'après le crit\`ere valuatif de propret\'e, il existe
un ouvert $O_1\subset Y_1$ 
dont le compl\'ementaire
est de codimension $\geq 2$ et une section
$\alpha_1$ au-dessus de $O_1$.

Le morphisme $Y_1\times_X Y_2 \ra Y_1$ est fini localement libre,
défini par l'image réciproque $\mathscr E_2$
par $Y_1\ra X$ de l'algèbre finie localement
libre qui définit $Y_2$.
La section $\alpha_1$ correspond ainsi à une section de $\mathscr E_2$
au-dessus de $O_1$. Comme $Y_1$ est \propS2 (le morphisme
$Y_1\ra S$ est \propS2 et plat comme composé de morphismes \propS2 et plats,
et $S$ est \propS2 puisque c'est un anneau de valuation discrète,
\cf EGA~4, prop.~6.8.3),
cette section s'étend en une unique section de $\mathscr E_2$ au-dessus
de $Y_1$, d'où une section du morphisme $Y_1\times_X Y_2\ra Y_1$,
d'où finalement un unique $X$-morphisme $\alpha:Y_1\ra Y_2$ qui
étend~$\alpha_\eta$.

En appliquant le même raisonnement à $\alpha_\eta^{-1}$,
puis à $\alpha_\eta\circ\alpha_\eta^{-1}$ et
$\alpha_\eta^{-1}\circ\alpha_\eta$,
on voit que $\alpha$ est un isomorphisme, comme il fallait démontrer.
\end{proof}

Il r\'esulte de la proposition pr\'ec\'edente et des
d\'efinitions le corollaire suivant :
\begin{coro}
Supposons que le morphisme $X\ra S$ satisfait à la propriété \propS2.
Le morphisme diagonal
$\mathcal H(X/S)\ra \mathcal H(X/S)\times_S \mathcal H$(X/S)
est représentable, fini et non ramifié.
\end{coro}

\begin{coro}
Supposons que le morphisme $X\ra S$ satisfait à la propriété \propS2.
Les champs $\mathcal H^n(X/S)$, $\mathcal H^{n,\text{$Z$-triv}}(X/S)$,
$\mathcal H^{n,\text{$U$-ét}}(X/S)$, $\mathcal H^G(X/S)$
sont des champs de Deligne--Mumford de type fini.
Le champ $\mathcal H(X/S)$ est un champ de Deligne--Mumford localement
de type fini.
\end{coro}
\begin{proof}
Rappelons qu'un champ de Deligne--Mumford
est un champ algébrique $\mathcal S$ qui possède une présentation
étale, c'est-à-dire tel qu'il existe un $1$-morphisme
$P\ra\mathcal S$, étale, où $P$ est un schéma.

C'est une conséquence de la prop.~\ref{prop.champ.alg.tf},
du corollaire précédent et de~\cite{laumon-mb2000}, th.~8.1.
\end{proof}

\begin{rema}
La théorie précédente s'étend \emph{verbatim} au cas d'un 
espace algébrique $X/S$ en remplaçant partout les schémas par des espaces
algébriques. Les seuls problèmes
proviennent de l'utilisation éventuelle de schémas de Hilbert 
(qui n'est un schéma que dans le cas projectif).
Cependant, étant donné un espace algébrique propre,
Artin construit dans~\cite{artin69}  son {\og espace algébrique
de Hilbert\fg}. La même démonstration que dans le cas projectif
implique alors que les fibrés vectoriels sur un espace algébrique propre
sont classifiés par un champ algébrique.
\end{rema}

\begin{rema}
Il faudrait aussi traiter un cas supplémentaire : celui où on fixe
des conditions sur les groupes d'inertie.
Ce sera essentiel au paragraphe suivant lorsqu'on voudra
exprimer les propriétés des champs de Hurwitz qui
résultent de SGA~1.
\end{rema}

\section{Propri\'et\'es des champs de Hurwitz}

Dans tout ce paragraphe, nous supposons que $X$ est lisse sur $S$
et que $U\subset X$ est l'ouvert complémentaire d'un diviseur
à croisements normaux relatifs $D$. Cela signifie (\cite{sga1}, XIII, 2.1)
que localement pour la topologie étale, $D\subset X$ est
isomorphe à $\div(x_1\cdots x_r) \subset \A^n$.

Nous allons donner quelques propriétés des champs de Hurwitz
en nous ramenant aux des énoncés correspondants sur les revêtements
ou sur les groupes fondamentaux démontrés dans~\cite{sga1}.
Si $f:(X,x)\ra (Y,y)$ est un morphisme de schémas connexes 
{\og géométriquement pointés\fg}, rappelons
que l'on a un homomorphisme $f_*:\pi_1(X,x)\ra \pi_1(Y,y)$
entre groupes fondamentaux et que :
\begin{itemize}
\item $f_*$ est injectif si et seulement si tout revêtement étale
connexe de $X$ est une composante connexe de l'image réciproque
par $f$ d'un revêtement étale de $Y$ (\cite{sga1}, V, 6.8) ;
\item $f_*$ est surjectif si et seulement si l'image réciproque
par $f$ de tout revêtement étale connexe de $Y$ est connexe 
(\cite{sga1}, V, 6.9) ;
\item $f_*$ est un isomorphisme si et seulement si $f$ induit par
image réciproque une bijection entre les catégories
des revêtements étales connexes
de $Y$ et de $X$ (conséquence des deux points précédents, \cite{sga1}, V,
6.10).
\end{itemize}

\begin{prop}
Le champ $\mathcal H^{\text{$U$-ét}}(X/S)$
est étale sur $S$ aux points correspondant
à des revêtements modérément ramifiés.
\end{prop}
\begin{proof}
Comme $\mathcal H$ est localement de type fini sur $S$, il suffit
de prouver qu'il est formellement étale aux points considérés.
Autrement dit, si $S'\subset S$ un fermé défini par un idéal de
carré nul et $Y'\ra X_{S'}$ un revêtement de $X_{S'}$ non ramifié
hors de $D_{S'}$ et modérément ramifié le long de $D_{S'}$,
on doit prouver qu'il existe
un unique revêtement $Y\ra X$ non ramifié hors de $D$
tel que $Y_{S'}=Y'$.

Or, c'est une conséquence du \emph{lemme d'Abhyankar relatif} (\cite{sga1},
XIII, 5.5), de l'énoncé analogue pour les revêtements étales
(\cite{sga1}, ??)
et de la théorie de la descente.
\end{proof}

Si $N$ est un entier, $S[1/N]$ désigne le plus grand ouvert de~$S$
où $N$ est inversible.

\begin{prop}\label{prop.finitude}
Le champ $\mathcal H^{G,\text{$U$-ét}}(X/S)$
est propre et quasi-fini au dessus de $S[1/\#G]$.
Le champ $\mathcal H^{n,\text{$U$-ét}}(X/S)\times_S S[1/n!]$
est propre et quasi-fini au-dessus de $S[1/n!]$.
\end{prop}
\begin{proof}
Posons $N\#G$ dans le premier cas et $N=n!$ dans le second.
Une fois établie l'assertion relative à la propreté, l'autre
découle de la proposition précédente : un $S$-champ de Deligne-Mumford
qui est propre et étale est en effet quasi-fini, comme on le voit
sur une présentation étale.

Le critère valuatif de propreté pour les champs
de Deligne-Mumford (\cite{laumon-mb2000}, prop.~7.12\footnote{La condition (*)
de \loccit est vérifiée car les champs considérés sont de Deligne-Mumford,
cf.~\loccit, remarques~7.12.2 et 7.12.3.})
nous ramène à établir le fait suivant :
supposons que $S=\Spec R$ est le spectre
d'un anneau de valuation discrète,
de corps résiduel $k$ de caractéristique $p$ ne divisant pas $N$ et de
corps des fractions $K$, alors tout revêtement galoisien
de groupe~$G$ (\resp de degré~$n$) de $X_K$ se prolonge
en un unique revêtement  de $X$ qui est galoisien de groupe~$G$
(\resp de degré~$n$).

Dans le premier cas, cela résulte de ce que le morphisme canonique
entre les plus grands quotients d'ordres premier à $p$
des groupes fondamentaux $\pi_1(X_K)^{(p)}\ra\pi_1(X)^{(p)}$
est injectif. Sa bijectivité résulte en fait de~\cite{sga1}, XIII, 2.8
(cf.~\emph{loc. cit.}, preuve de 2.12, page 394).

Dans le second cas, l'ordre du groupe de Galois de la clôture galoisienne
d'un tel revêtement divise $n!$, ce qui nous ramène au cas déjà traité.
\end{proof}

\begin{coro}\label{coro.schemas}
Les champs de Hurwitz $\mathcal H^{G,\text{$U$-ét}}(X/S)$,
$\mathcal H^{n,\text{$U$-ét}}$
que nous avons introduits possèdent des espaces grossiers en tant
qu'espaces algébriques.
Leur restriction au plus grand ouvert tel que $\# G$ (\resp $n!$)
est inversible est un schéma fini étale.
\end{coro}
\begin{proof}
D'après un théorème de Keel et Mori (\cite{keel-m1997}, cor.~1.3;
voir aussi~\cite{laumon-mb2000}, th.~19.1),
tout champ algébrique de type fini
et séparé possède un espace grossier associé dans la
catégorie des espaces algébriques. Notons ainsi avec des
lettres romaines les $S$-espaces algébriques grossiers
des champs de Hurwitz.
Il résulte de la proposition~\ref{prop.finitude}, 
que $\mathrm H^{G,\text{$U$-ét}}(X/S)$ (\resp $\mathrm
H^{n,\text{$U$-ét}}(X/S)$) est propre et quasi-fini au-dessus de
$S[1/N]$, où $N=\# G$ (\resp $N=n!$).
La proposition est alors une conséquence immédiate de ce
qu'un espace algébrique quasi-fini au-dessus d'un schéma
est un schéma (\cite{laumon-mb2000}, th.~A.2 et cor.~A.2.1).
\end{proof}

Bien sûr, dans les cas où le problème de modules
considéré n'a pas d'automorphisme non trivial,
la proposition précédente montre que le champ de Hurwitz
est est un schéma fini étale au-dessus d'un ouvert convenable
de~$S$.
C'est notamment le cas du champ $\mathscr H^G(X/S)$,
lorsque le centre du groupe fini~$G$ est trivial.

\section{\`A propos d'un théorème de S.~Beckmann}

\paragraph{}
Soit $S$ un schéma et soit $\mathscr X$ un $S$-champ algébrique
localement noethérien.
Rappelons (\cite{laumon-mb2000}, chap.~5 et 11) qu'un point
(automatiquement {\og algébrique\fg}, en vertu de \loccit, th.~11.3)
de~$X$ est une classe d'équivalence de couples $(x,K)$, où
$K$ est un $S$-corps\footnote{Un $S$-corps est un corps~$K$ 
mundi d'un morphisme de~$\Spec K$ vers~$S$.}
et $x$ un morphisme de $S$-champs, $\Spec K\ra \mathscr X$.
On dit que deux couples $(x_1,K_1)$ et $(x_2,K_2)$ sont équivalents
s'il existe un couple $(x,K)$, où $K$ est un sur-$S$-corps de $K_1$
et $K_2$ et $x=x_1\otimes_{K_1} K=x_2\otimes_{K_2} K$.

Un tel point~$\xi$ se factorise canoniquement à travers un sous-$S$-champ
$\mathscr G_{\xi}$ de~$\mathscr X$ qui est
une gerbe de type fini sur son faisceau grossier, lequel est
de la forme $\Spec \kappa(\xi)$, pour un certain corps~$\kappa(\xi)$
appelé \emph{corps résiduel du point~$\xi$}.

On voit ainsi que, dans le langage des espaces de modules,
le corps résiduel d'un point n'est autre que le \emph{corps des modules}
de l'objet représenté.

Nous pouvons maintenant énoncer l'analogue pour les champs
algébriques d'un théorème de Chevalley et Weil :

\begin{prop}\label{prop.c-w}
Soit $S$ un schéma normal, intègre, de corps des fractions $K$.
Soit $\mathscr X$ un $S$-champ de Deligne-Mumford
qui est propre et étale sur~$S$.
Alors, le corps résiduel de tout point algébrique du $K$-champ $\mathscr X_K$
est non ramifié sur~$S$.
\end{prop}
\begin{proof}
Pour la démonstration,
on peut supposer que $S$ est le spectre d'un anneau de valuation
discrète.

Soit $X$ l'espace algébrique grossier associé à $\mathscr X$.
D'après les arguments du cor.~\ref{coro.schemas}, c'est en fait un
$S$-schéma.
La factorisation $\mathscr X\ra X\ra S$ et le fait que $\mathscr X$
soit propre et étale impliquent que $X$ est fini et étale sur~$S$.

Soit $\xi$ un point de~$\mathscr X$. Son image dans~$X$ définit
un $\kappa(\xi)$-point~$x$ d'un $S$-schéma fini et étale, où le
corps résiduel de~$\xi$, $\kappa(\xi)$, 
est une extension finie de~$K$. Soit $R$ la clôture intégrale de~$S$
dans~$\kappa(\xi)$. D'après le critère valuatif de propreté,
$x$ se prolonge en une section $\Spec R\ra X$ dont l'image
est nécessairement étale sur~$S$. Autrement dit,
$\kappa(\xi)$ n'est pas ramifiée sur~$K$.
\end{proof}

Une première conséquence de ces considérations est une généralisation
d'un théorème de S.~Beckmann~\cite{beckmann89}
concernant le corps des modules d'un
revêtement ramifié de la droite projective.

\begin{theo}
Soit $S$ un schéma normal, intègre, de corps des fractions~$K$;
soit $X$ un $S$-schéma lisse, $D$ un diviseur à croisements normaux
relatif dans~$X$. 
Soit $G$ un groupe fini (\resp soit $n$ un entier); on suppose
que le cardinal de~$G$ (\resp $n!$) est inversible sur~$S$.

Alors, le corps des modules d'un revêtement (a priori
défini sur une extension finie de~$K$) galoisien de groupe~$G$
(\resp d'un revêtement de degré~$n$) de~$X$, non ramifié hors de~$D$,
est non ramifié sur~$S$.
\end{theo}

\begin{proof}
Il suffit d'appliquer le théorème de Chevalley-Weil (prop.~\ref{prop.c-w})
au champ de Hurwitz adéquat ($\mathscr H^{G, \text{$X\setminus
D$-ét}}(X/S)$, \resp $\mathscr H^{n, \text{$X\setminus D$-ét}}(X/S)$)
et au point fourni par le revêtement.
\end{proof}

Notons que dans un autre article (\cite{beckmann91}), Beckmann 
a établi un théorème qui concerne cette fois la ramification
de l'extension de corps de nombres fournie par la fibre d'un tel revêtement.
Une version plus faible découle du théorème de Chevalley-Weil
appliqué au revêtement lui-même. Notons aussi que Conrad
a généraliser ce second théorème dans~\cite{conrad2000}
au cas de la dimension supérieure. Il serait probablement intéressant
d'interpréter ce résultat dans le cadre de la géométrie logarithmique.


\providecommand{\noopsort}[1]{}\providecommand{\url}[1]{\textit{#1}}
\providecommand{\bysame}{\leavevmode ---\ }
\providecommand{\og}{``}
\providecommand{\fg}{''}
\providecommand{\smfandname}{\&}
\providecommand{\smfedsname}{\'eds.}
\providecommand{\smfedname}{\'ed.}
\providecommand{\smfmastersthesisname}{M\'emoire}
\providecommand{\smfphdthesisname}{Th\`ese}

\end{document}